\documentstyle[12pt]{article}
\newlength{\spacing}
\newcommand{\doublespace}{\setlength{\baselineskip}{1.5\spacing}}
\newtheorem{thm}{Theorem}[section]

\newtheorem{defn}{Definition}

\def\rar{\to}

\def\del{\delta}

\def\ep{\epsilon}

\def\today{\ifcase\month\or
  January\or February\or March\or April\or May\or June\or
  July\or August\or September\or October\or November\or December\fi
  \space\number\day, \number\year}

\begin{document}
\begin{titlepage}
\begin{center}
%
%
{\bf A Strong threshold for the size of random caps to cover a sphere.} \\
\vspace{0.20in} by \\
\vspace{0.2in} {Bhupendra Gupta \footnote{Corresponding Author.
email:gupta.bhupendra@gmail.com, bhupen@iiitdm.in}}\\
Faculty of Engineering and Sciences,\\
Indian Institute of Information Technology (DM)-Jabalpur, India.\\
\vspace{0.1in}
\end{center}
\vspace{0.2in}
\sloppy
\begin{center} {\bf Abstract} \end{center}

\begin{center} \parbox{4.8in}
{In this article, we consider `$N$'spherical caps of area $4\pi p$
were uniformly distributed over the surface of a unit sphere. We are
giving the strong threshold function for the size of random caps to
cover the surface of a unit sphere. We have shown that for large
$N,$ if $\frac{Np}{\log\:N} > 1/2$ the surface of sphere is
completely covered by the $N$ caps almost surely , and if
$\frac{Np}{\log\:N} \leq 1/2$ a partition of the surface of sphere
is remains uncovered by the $N$ caps almost surely. }
\\
\vspace{0.4in}
\end{center}

\vspace{0.5in} {\sl AMS 2000 subject classifications}:
\hspace*{0.5in} 05C80, 91D30.
\end{titlepage}
\doublespace
\section{Introduction.}
Let $V_1,V_2,\ldots,V_N$ be the spherical caps on the surface of a
unit sphere with their centers $v_1,v_2,\ldots,v_N$ respectively, on
the surface of a unit sphere. Also let $v_1,v_2,\ldots,v_N$ are
independently and uniformly distributed on the surface of a unit
sphere. H. Maehara, \cite{maehara} gives the threshold function
$p_0(N) = \frac{\log\:N}{N}$ for the coverage of the surface of a
unit sphere. H. Maehara, proves that for $\frac{p(N)\cdot
N}{\log\:N}< 1,$ probability that $N$ spherical caps cover the
entire surface of the unit sphere is converges to $0,$ and for
$\frac{p(N)\cdot N}{\log\:N}>1,$ probability that each point of the
sphere is covered by $n$ caps is converges to $1.$ Since both of
these are convergence in probability sense, the threshold $p_0(N)$
is a weak threshold. Also in article \cite{maehara}, instead of
exact bounds author use lose approximations. Due to these
approximations threshold suggested in \cite{maehara} is different
from threshold suggested in this article.\\

Now using the same model and notations as in H. Maehara,
\cite{maehara}, we are giving the strong threshold function
for the coverage of the surface of a unit sphere.\\

\section{Basic Model and Definitions. }

Here, we recall the same model as it is given by H. Maehara,
\cite{maehara}. We made some modification in the language for making
the things more clear.

Let `$S$' be the surface area of a unit sphere in $3-$dimensional
space. Let $V_1,V_2,\ldots,$ be the spherical caps on the surface of
a unit sphere with their centers $v_1,v_2,\ldots,$ respectively, and
uniformly distributed on the surface of a unit sphere. The area of a
spherical cap of angular distance (angular radius) `$a$' is $2\pi
(1-\cos(a)) = 4\pi \sin^2(a/2).$

Let `$p$' be the probability that any point on the surface of unit
sphere covered by a specified spherical cap of angular distance
(angular radius) `$a$'. Then
\begin{equation}
p :=  \frac{\mbox{Area of a spherical cap of angular distance `$a$'}
}{\mbox{Surface area of unit sphere}} = \frac{4\pi
\sin^2(a/2)}{4\pi} = \sin^2(a/2)\label{p}
\end{equation}
Let there are `$N$' random caps of angular distance `$a$' on the
surface `$S$' of unit sphere. Let $U_0(N,p),$ be the set of those
points which remains uncovered by `$N$' spherical caps and
$u_0(N,p)$ be the proportion of the area covered by $U_0(N,p),$
i.e., $u_0(N,p)$ be the proportion of the area which remains
uncovered by `$N$' spherical caps:
\begin{equation}
u_0(N,p) := \frac{\{\mbox{the area of }
U_0(N,p)\}}{4\pi}\:.\label{u_n}
\end{equation}
Then
\begin{equation}
E(u_0(N,p)) = \frac{1}{4\pi}\int_SP[x \in U_0(N,p)]dx.\label{Eu}
\end{equation}
Now consider,
\begin{eqnarray}
P[x \in U_0(N,p)] & = & P[x \mbox{ is remains uncovered }]\nonumber\\
& = & \prod_{i=1}^{N}\left(1-P[x \in V_i, ]\right)\nonumber\\
& = & (1-p)^N.\label{x_in_U}
\end{eqnarray}
Hence, from (\ref{Eu}), we have
\begin{equation}
 E(u_0(N,p)) = (1-p)^N. \label{Eu_n}
\end{equation}
Similarly, as in H. Maehara, \cite{maehara}, we have
\begin{equation}
E(u^2_0(N,p)) = \frac{1}{16\pi^2}\int_S\int_SP[x,y \in U_0(N,p)]dx =
\frac{1}{4\pi}\int_SP[x_0,y \in U_0(N,p)]dy, \label{Eu_n2}
\end{equation}
where $x_0$ is a fixed point on $S.$ Let $x_0$ and $y$ subtend an
angle `$\theta$' at the center of sphere. Then
\begin{equation}
P[x_0,y \in U_0(N,p)] = \left(1-(2p-q(\theta))\right)^N,
\end{equation}
where, $q(\theta)$ be the area of intersection between two spherical
caps of angular distance `$a$'. Substituting the above probability
in (\ref{Eu_n2}), we get
\[E(u^2_0(N,p)) = \frac{1}{4\pi}\int_S\left(1-(2p-q(\theta))\right)^Ndy.\]
Since points $x_0$ and $y$ subtend an angle between $\theta$ and
$\theta +d\theta$ at the center of the sphere. Then
\begin{equation}
E(u^2_0(N,p)) = \int_{0}^{\pi}\left(1-(2p-q(\theta))\right)^N
(1/2)\sin(\theta)d\theta.\label{eu2}
\end{equation}
Since $q(\theta) = 0$ for $\theta >2a,$
\begin{eqnarray}
E(u^2_0(N,p)) & < & \int_{0}^{2a}(1-p)^N(1/2)\sin(\theta)d\theta
+\int_{2a}^{\pi}(1-2p)^N(1/2)\sin(\theta)d\theta\nonumber\\
& < & (1-p)^N[-(1/2)\cos(\theta)]_{0}^{2a} + (1-2p)^N\nonumber\\
& = & (1-p)^N\frac{1-\cos(2a)}{2} + (1-2p)^N.\label{x1}
\end{eqnarray}
Using (\ref{Eu_n}), in (\ref{x1}), we have
\begin{eqnarray}
E(u^2_0(N,p)) & < & (1-p)^N\frac{1-\cos(2a)}{2} +
(1-2p)^{N}\nonumber\\
& = & 4p(1-p)^{N+1}+(1-2p)^{N},\label{eu2uper}
\end{eqnarray}
since $\frac{1-\cos(2a)}{2} = 4p(1-p).$\\

Now for the lower bound of $E(u^2_0(N,p)),$ from (\ref{eu2}) we have
\begin{equation}
E(u^2_0(N,p))  > \int_{2a}^{\pi}\left(1-(2p-q(\theta))\right)^N
(1/2)\sin(\theta)d\theta,
\end{equation}
since $q(\theta) = 0$ for $\theta >2a.$
\begin{eqnarray}
E(u^2_0(N,p))
& > & \frac{(1-2p)^N}{2}\int_{2a}^{\pi}\sin(\theta)d\theta\nonumber\\
& = & \frac{(1-2p)^N}{2}[-\cos(\theta)]_{2a}^{\pi}\nonumber\\
& = & (1-2p)^N\frac{\cos(2a)+1}{2}\nonumber\\
& = & (1-2p)^N(1-4p(1-p)),\label{eu2lower}
\end{eqnarray}
since $\frac{\cos(2a)+1}{2} = 1-4p(1-p).$\\

Let $\Theta$ be a fixed monotone property.

\begin{defn}

  A function $\del_{\Theta}(c): Z^+ \rar {R}^+$ is a {\it strong
  threshold}
  function for $\Theta$ if the following is true for every fixed
  $\ep>0,$
 \begin{itemize}
  \item if $P[\del_{\Theta}(c-\ep) \in \Theta] = 1 - o(1),$ and
  \item if $P[\del_{\Theta}(c+\ep) \in \Theta] = o(1)\ ,$
  \end{itemize}
where `$c$' is some constant.  \hfill$\Box$\\
\end{defn}
\section{Main Result.}
\begin{thm}
Let $p = \frac{c\log\:N}{N},$ then for $c> \frac{1}{2},$ we have the
surface of unit sphere is completely covered by `$N$' spherical
caps, i.e.,
\[U_0(N,p) = \phi,\qquad \mbox{almost surely},\]
and for $c\leq \frac{1}{2},$ surface of unit sphere is not
completely covered by `$N$' spherical caps, i.e.,
\[U_0(N,p) \neq \phi,\qquad \mbox{almost surely.}\]
\end{thm}
\textbf{Proof.} For arbitrary small $\ep.$ We have,
\begin{equation}
P[U_0(N,p) \neq \phi] \simeq P[\mid u_0(N,p)\mid \geq \ep].
\label{e3}
\end{equation}
From the Markov's inequality, we have
\begin{eqnarray}
%
P[ \mid u_0(N,p)\mid \geq \ep] & \leq &
\frac{E[u^2_0(N,p)]}{\ep^2}\nonumber\\
& < & \frac{1}{\ep^2}\left(4p(1-p)^{N+1}+(1-2p)^{N}\right),
\end{eqnarray}
using the upper bound of $E[u^2_0(N,p)]$ from (\ref{eu2uper}). Now
taking $p = \frac{c\log\:N}{N},$ where $c$ is some constant. Then
\begin{eqnarray}
P[ \mid u_0(N,p)\mid \geq \ep] & < &
\frac{1}{\ep^2}\left(\frac{4c\log\:N}{N}\left(1-\frac{c\log\:N}{N}\right)^{N+1}+
\left(1-\frac{2c\log\:N}{N}\right)^N\right)\nonumber\\
& < &
\frac{1}{\ep^2}\left(\frac{4c\log\:N}{N}\left(1-\frac{c\log\:N}{N}\right)e^{-c\log\:N}+
e^{-2c\log\:N}\right)\nonumber\\
& < &
\frac{1}{\ep^2}\left(\frac{4c\log\:N}{N^{1+c}}+\frac{1}{N^{2c}}\right).\label{e1}
\end{eqnarray}
If $c > 1/2,$ the above probability is summable, i.e.,
\[\sum_{N=0}^{\infty} P[\mid u_0(N,p)\mid \geq \ep] < \infty,\]
and hence from (\ref{e3}), we have
\[\sum_{N=0}^{\infty}P[U_0(N,p) \neq \phi] < \infty.\]
Then by the Borel-Cantelli's Lemma, we have
%
\[P[U_0(N,p) \neq \phi, \qquad i.o.] = 0.\]
Thus, the set $U_0(N,p)\neq \phi$ for only finitely many time, i.e.,
eventually $U_0(N,p)= \phi$ happens infinitely times with
probability $1.$
Hence for $c > \frac{1}{2},$ we have
\[U_0(N,p) = \phi,\qquad \mbox{ almost surely.} \]\\

Now, by the lower bound of Chebyshev's inequality (Page 55,
Shiryayev \cite{Shiryayev}.), we have
\begin{equation}
P[ u_0(N,p) \geq \ep] \geq \frac{E[ u_0(N,p)^2]-\ep^2}{16\pi^2},
\end{equation}
since $u_0(N,p)\geq 0$ and $\mid u_0(N,p) \mid \leq 4 \pi.$ Now
using the lower bound of $E[ u_0(N,p)^2]$ from (\ref{eu2lower}) and
taking $\ep = o\left(\frac{1}{N}\right),$ we have
\begin{eqnarray}
P[ u_0(N,p) \geq \ep] & \geq &
\frac{(1-2p)^N(1-4p(1-p))-\ep^2}{16\pi^2}\nonumber\\
& \geq & C_1(1-2p)^N(1-4p(1-p))-C_2,
\end{eqnarray}
where $C_1= \frac{1}{16\pi^2}$ and $C_2= \frac{\ep^2}{16\pi^2}.$\\

Substituting $p = \frac{c\log\:N}{N},$ in the above expression we
get
\begin{eqnarray}
P[u_0(N,p) \geq \ep] & \geq &
C_1\left(1-\frac{2c\log\:N}{N}\right)^N\left(1-\frac{4c\log\:N}{N}\left(1-\frac{c\log\:N}{N}\right)\right)-C_2\nonumber\\
& \geq & C_3e^{-2c\log\:N} = \frac{C_3}{N^{2c}},\label{e2}
\end{eqnarray}
where $C_3$ is some constant. If we take $c \leq 1/2,$ then the
probability (\ref{e2}), is not summable with respect to $N,$ i.e.,
\[\sum_{N=1}^{\infty}P[ u_0(N,p) \geq \ep] = \infty.\]
Then by the Borel-Cantelli's Lemma, we have
\[P[u_0(N,p) \geq \ep, \qquad i.o.] = 1,\]
since $u_0(N,p)$ are independent. Thus $u_0(N,p) \geq \ep$ happens
infinitely many time with probability $1.$ Hence for $c \leq
\frac{1}{2},$ we have
\[  u_0(N,p) \geq \ep, \qquad \mbox{almost surely.}\]
This implies for $c\leq 1/2,$ we have $U_0(N,p) \neq \phi$ almost surely.\hfill$\Box$\\

\end{document}